\documentclass{article}

\newtheorem{fed}{\textbf{Definition}}[section]
\newtheorem{thm}[fed]{\textbf{Theorem}}
\newtheorem{lemma}[fed]{\textbf{Lemma}}

\newtheorem{rem}[fed]{\textbf{Remark}}
\newtheorem{prop}[fed]{\textbf{Proposition}}
\newtheorem{cor}[fed]{\textbf{Corollary}}
\usepackage{amssymb,bbm,graphicx,epsfig,psfrag,epic,eepic,latexsym}
\usepackage{amsmath}
\usepackage{mathrsfs} 
\usepackage[dvips]{color}       

\begin{document}
\title{Finite dimensional approximations for the symplectic vortex equations}
\author{Urs Frauenfelder\footnote{
Supported by Swiss National Science
Foundation and JSPS}}
\maketitle

\begin{abstract}
In this paper we study Furuta's finite dimensional approximations for the flow
lines of the Moment action functional for toric symplectic orbifolds.
These are the analogon of the finite dimensional approximations for the
Chern-Simons-Dirac functional studied by Kronheimer and Manolescu, see
\cite{kronheimer-manolescu, manolescu}.
The Moment action functional is the sum of Floer's action functional and
a Lagrangian multiplier to a constraint given by the moment map. Its
flow lines are the symplectic vortex equations on the cylinder.   

We prove a compactness theorem specific for functions containing a
Lagrange multiplier which allows us to define the Conley indices of the
flow on the finite dimensional approximations. We show that these Conley
are given by the Thom space of the normal bundle of toric map space
introduced by Givental in 
\cite{givental2}. Finally, we show how the moduli
spaces of the Morse flow lines on the finite dimensional approximations
are related to the moduli spaces of the symplectic vortex equations on
the cylinder.   
\end{abstract}
\tableofcontents

\newpage

\section[Introduction]{Introduction}

In \cite{kronheimer-manolescu, manolescu} finite dimensional
approximation for the Seiberg Witten equation were studied. In this
paper we will study their analogon for the symplectic vortex
equations of toric symplectic orbifolds. The symplectic vortex equations
were introduced in \cite{cieliebak-gaio-mundet-salamon, 
cieliebak-gaio-salamon,mundet} and their solutions appeared in
\cite{frauenfelder1,frauenfelder2} as flow lines in moment Floer 
homology.
We prove that the finite dimensional approximations exist and show their
relation to a finite dimensional model for the moduli spaces
of Floer homology for toric symplectic
orbifolds considered before by Givental in \cite{givental2}, see also
\cite{givental, iritani,vlassopoulos}. To examine the relation between
the original equations and their finite dimensional approximations we
prove a cobordism theorem which relates the moduli spaces of the
symplectic vortex equations to the moduli spaces of finite dimensional
Morse theories.

The idea behind these finite dimensional approximations is the
following. Assume that we have an action functional
$\mathcal{A}$ and a metric 
defined on an infinite dimensional space $\mathscr{L}$.
Suppose further that there are finite dimensional spaces $L_\nu$ which
approximate $\mathscr{L}$. We consider the sets $K_\nu \subset L_\nu$
which are given by the traces of all the flow lines of the action
functional restricted to $L_\nu$ which converge at both ends, i.e.
$$K_\nu=\{y(\sigma): \sigma \in \mathbb{R},\,\,
y \in C^\infty(\mathbb{R},L_\nu),\,\,
\partial_s y+\nabla \mathcal{A}|_{L_\nu}=0,\,\,
\exists \lim_{s\to \pm \infty}y(s)\}.$$
The sets $K_\nu$ are clearly invariant under the local flow of the restricted
action functional. If they 
are also compact one can associate to them a
Conley index $I_\nu$. For good approximations one expects that the union
of these Conley indices should contain all the information of the moduli
space of all finite energy gradient trajectories of the action
functional $\mathcal{A}$. As it was explained in
\cite{kronheimer-manolescu} these Conley indices should be the building
blocks of Furuta's pro-spectrum with parametrized universe which in the
unobstructed case reduces to the Floer pro-spectrum introduced by Cohen,
Jones and Segal in \cite{cohen-jones-segal}.

In \cite{kronheimer-manolescu, manolescu} Kronheimer and Manolescu
studied finite dimensional approximations for the Chern-Simons-Dirac 
action functional which were used before by Furuta in his celebrated 
proof of the 10/8-conjecture, see \cite{furuta}.
In \cite{manolescu} Manolescu studied the case where the
Chern-Simons-Dirac action functional was bounded on its
critical set. Using this he could get compactness for the
sets $K_\nu$. In \cite{kronheimer-manolescu} Kronheimer and Manolescu
considered the case where the Chern-Simons-Dirac action functional was
not any more bounded on its critical set. In this case it is not clear
if the sets $K_\nu$ are compact. 
However, Kronheimer and Manolescu could still
define something provided a twist discovered by Furuta vanishes. In this
paper we consider Furutas finite dimensional approximations for the
action functional of moment Floer homology. This action functional will
in general not be bounded on its critical set. However, we will prove
that we have compactness on the finite dimensional approximations.

Assume that a Lie group acts on a symplectic manifold by Hamiltonian 
isomorphisms, i.e. there exists a moment map $\mu$ for the action. Then
the Moment action functional $\mathcal{A}$
is defined as Floer's
action functional together with a Lagrangian multiplier to the
constraint $\mu^{-1}(\tau)$, where $\tau$ is an element of the Lie
algebra of the acting Lie group. In this paper we consider as
symplectic manifold the complex vector space $\mathbb{C}^n$ on which the
torus $T^k$ acts via a representation to the unitary group. 
We put the following standing assumption
\begin{description}
 \item[(H)] The moment map $\mu$ is proper and $\tau \in \mathfrak{t}^k$ is a
  regular value of $\mu$.
\end{description}
Under this hypothesis the 
Marsden-Weinstein quotient $\mu^{-1}(\tau)/T^k$ 
is a compact toric symplectic orbifold. 

We denote by $\mathscr{L}$ the space of smooth loops on
$\mathbb{C}^n  \times \mathfrak{t}^k$ where $\mathfrak{t}^k$ is the
Lie algebra of the torus $T^k$. Then 
$$\mathcal{A}\colon \mathscr{L} \to \mathbb{R}$$
and there is an action on $\mathcal{L}$
of the gauge group $\mathcal{H}$ consisting of smooth loops on the torus under
which the diffential of the action functional is invariant. There is a
natural splitting
$$\mathcal{H}=H \times \mathcal{H}_0$$
such that the group $H \cong T^k \times \mathbb{Z}^k$ is finite, the group
$\mathcal{H}_0$ is contractible and acts freely on $\mathscr{L}$, 
and the action functional
is invariant under $\mathcal{H}_0$. Since the gauge group is abelian we
can find a global Coulomb section $\mathscr{L}_0$ in the principal
$\mathcal{H}_0$-bundle $\mathscr{L}$. This section is $H$-invariant, i.e.
the following diagram commutes and is equivariant under the group $H$.

\setlength{\unitlength}{1.5cm}
\begin{picture}(5,2)\thicklines
 \put(2,0){$\mathscr{L}_0$}
 \put(2.3,0.3){\vector(2,1){2}}
 \put(3.2,0.9){$\iota$}
 \put(2.5,0.1){\vector(1,0){1.7}}
 \put(3.4,0,2){$c$}
 \put(4.3,0){$\mathscr{L}/\mathcal{H}_0$}
 \put(4.4,1.2){$\mathscr{L}$}
 \put(4.5,1.1){\vector(0,-1){0.8}}
 \put(4.6,0.6){$\pi$}
\end{picture}
\\ \\
The $L^2$-inner product endows $\mathscr{L}$ with an $\mathcal{H}$-invariant
metric $g$. We denote by $\bar{g}$ its quotient metric on 
$\mathscr{L}/\mathcal{H}_0$. There are two natural $H$-invariant metrics on
$\mathscr{L}_0$
$$g_0=\iota^* g, \quad g_1=c^* \bar{g}.$$ 
We will abbreviate
$$\mathcal{A}_0=\iota^*\mathcal{A}.$$
Note that flow lines of the gradient of $\mathcal{A}$ with respect to
$g$ are in one to one correspondence with flow lines of the gradient of 
$\mathcal{A}_0$ with respect to $g_1$.

We approximate $\mathscr{L}_0$ by finite dimensional submanifolds
$L_\nu$. These are of the form
$$L_\nu \cong \mathbb{C}^{N_\nu}\times \mathfrak{t}^k.$$
The restriction of the action functional is linear in $\mathfrak{t}^k$,
i.e. there exist smooth functions 
$f_\nu, h_\nu \in C^\infty(\mathbb{C}^{N_\nu},\mathbb{R})$ such that
$$\mathcal{A}_0|_{L_\nu}(z,\eta)=f_\nu(z)+\langle h_\nu(z),\eta\rangle,
\quad z \in \mathbb{C}^{N_\nu},\,\,\eta \in \mathfrak{t}^k.$$
In particular, the critical points of $\mathcal{A}_0|_{L_\nu}$ are
critical points of $f_\nu$ to the Lagrangian constraint $h_\nu^{-1}(0)$. 
It turns out that the  quotient of the zero sets of $h_\nu$ 
under the torus action are precisely the
toric map spaces which were introduced
before by Givental in \cite{givental2}
as a finite dimensional approximation to the moduli spaces of Floer
homology for the symplectic orbifold $\mu^{-1}(\tau)/T^k$, see also
\cite{givental, iritani, vlassopoulos}.

We consider the finite energy gradient flow lines of 
$\mathcal{A}_0|_{L_\nu}$ with respect to the metrics $g_0$ and $g_1$. We
show that its traces are contained in a compact set which allows us to
define the corresponding Conley indices $I_\nu^{g_0}$ and $I_\nu^{g_1}$.
More precisely, our first main result is the following theorem.
\\ \\
\textbf{Theorem A} \emph{The Conley indices $I_\nu^{g_0}$ and
 $I_\nu^{g_1}$ are well defined and isomorphic to the Thom space of
the normal bundle of $h_\nu^{-1}(0)$ in $\mathbb{C}^{N_\nu}$. In
particular, they are isomorphic to each other.}
\\ \\
The idea of Theorem A is the following. For $r \in [0,1]$ we 
consider the family of functions 
$F_{\nu,r} \in C^\infty(\mathbb{C}^{N_\nu}\times
\mathfrak{t}^k,\mathbb{R})$
defined by
$$F_{\nu,r}(z,\eta)=r f_\nu(z)+\langle h_\nu(z),\eta\rangle.$$
Then $F_{\nu,1}=\mathcal{A}_0|_{L_\nu}$ and for $r=0$ all flow lines of
$F_{\nu,0}$ are constant and lie in $h_\nu^{-1}(0) \times \{0\}$.

We will denote by $\nabla _0 \mathcal{A}_0$ the gradient of 
$\mathcal{A}_0$ with respect to the metric $g_0$.
The submanifolds $L_\nu$ of $\mathscr{L}_0$ are totally geodesic with
respect to the metric $g_0$ and if $p \in L_\nu$ then
$$\nabla_0\mathcal{A}_0(p) \in T_p L_\nu.$$
Hence every flow line of the restriction of $\mathcal{A}_0$ to $L_\nu$
with respect to $g_0$ is actually a flow line of $\mathcal{A}_0$
on $\mathscr{L}$ with respect to $g_0$. We will show that each flow line
of finite energy of $\nabla_0\mathcal{A}_0$ is actually contained in
some finite dimensional approximation.
\\ \\
\textbf{Theorem B }\emph{Assume that $y \colon \mathbb{R} \to
\mathscr{L}_0$ satisfies 
$$\partial_s y(s)+\nabla_0\mathcal{A}_0(y(s))=0,\,\, s\in
\mathbb{R},\quad \exists \lim_{s \to \pm \infty}y(s).$$
Then there exists $L_\nu$ such that the trace of $y$ is contained
entirely in $L_\nu$, i.e. 
$$y(s) \in L_\nu, \,\,s \in \mathbb{R}.$$}
\\ 
Finally we prove a cobordism theorem between the moduli spaces of the 
gradient flow lines
of $\mathcal{A}_0$ with respect to the metric 
$g_0$ and the ones with 
respect to the metric $g_1$. More precisely, we show the following
theorem.
\\ \\
\textbf{Theorem C }\emph{The metrics $g_0$ and $g_1$ on $\mathscr{L}_0$
can be connected by a path of $H$-invariant
metrics $g_r$ for $r \in [0,1]$ with the
following property. Denote by $\nabla_r \mathcal{A}_0$ the gradient of
$\mathcal{A}_0$ with respect to the metric $g_r$. Assume that 
$y_\nu$ for $\nu \in \mathbb{N}$ is a sequence of solutions of 
$$\partial_s y_\nu(s)+\nabla_{r_\nu} \mathcal{A}_0(y(s))=0, \quad s
\in \mathbb{R},\,\,r_\nu \in [0,1]$$
whose energy is uniformly bounded. Then there exists a subsequence
$\nu_j$, a sequence of gauge transformations $h_j \in H$ and a flow line
$y$ of $\nabla_r \mathcal{A}_0$ for $r \in [0,1]$, such that
$(h_j)_*y_{\nu_j}$ converges with respect to the
$C^\infty_{loc}$-topology to $y$, i.e.
$$(h_j)_*y_{\nu_j}\stackrel{C^\infty_{loc}}{\longrightarrow}_{j \to
\infty} y.$$}
\\ \\
The metric $g_r$ for $r \in (0,1]$ will be constructed in the following
way. We deform the action of $\mathcal{H}_0$ on $\mathscr{L}$ in an
$H$-invariant way by free actions for which $\mathscr{L}_0$ is still a
global section. The metric $g_r$ is then defined by pulling back the
quotient metric on $\mathscr{L}/_r \mathcal{H}_0$ induced by the
original $L^2$-metric $g$ on $\mathscr{L}$.
\\

This paper is organized as follows. In Section~\ref{toric} we recall the 
definition of toric symplectic orbifolds and prove a regularity
criterion which will allow us later to prove the regularity of
Givental's toric map spaces. 

In Section~\ref{moment} we recall the definition of the Moment action 
functional. We define the global Coulomb section on the loop space and
describe the two metrics $g_0$ and $g_1$. We compute the gradient of the
action functional with respect to both metrics and prove an estimate for
the two metrics which will allow us later to conclude that the finite
dimensional approximations are geodesically complete.

In Section~\ref{finite} we introduce the finite dimensional
approximations. We prove a regularity result for Givental's toric map space
and show
that the finite dimensional approximations are geodesically complete for
both metrics. These properties will be crucial to prove that the Conley
index is well defined. 

In Section~\ref{conley} we prove Theorem A. We first introduce the
Conley indices. Then we will prove a more general theorem which allows
us to compute the Conley index for a large class of functions which
contain some Lagrange multiplier part. Theorem A will follows as an
easy Corollary of this more general theorem and the results obtained in
Section~\ref{finite}. 

In Section~\ref{cobordism} we will prove Theorems B and C. 
\\ \\
\textbf{Acknowledgements: }I would like to express my deep gratitude to
K.Ono for pointing my attention to the work of Furuta and
Kronheimer-Manolescu. I would like to thank him and H.Iritani for useful
discussions.

\section[Toric symplectic manifolds]{Toric symplectic orbifolds}\label{toric}

Assume that for $k\leq n$ the torus 
$T^k=\{e^{iv}: v \in \mathbb{R}^k\}$ 
acts on the complex vector space
$\mathbb{C}^n$ via the action
$$\rho(e^{iv})z=e^{iAv}z, \quad z \in \mathbb{C}^n, \,\,v \in \mathbb{R}^k$$
for some $(n\times k)$-matrix $A$ with integer entries. We endow the 
Lie algebra of the torus
$$\mathrm{Lie}(T^k)=\mathfrak{t}^k=i\mathbb{R}^k$$
with its standard inner product. 
The action of the torus
on $\mathbb{C}^n$ is Hamiltonian with respect to the standard symplectic
structure $\omega=\sum_{i=1}^n dx_i \wedge dy_i$. 
Denoting by $A^T$ the transposed matrix of $A$ a moment map 
$\mu \colon \mathbb{C}^n \to \mathfrak{t}^k$ is given by
\begin{equation}\label{moment}
\mu(z)=\frac{i}{2}A^Tw, \quad w=\left(\begin{array}{c}
|z_1|^2 \\
\vdots\\
|z_n|^2
\end{array}\right),
\end{equation}
i.e.
$$d\langle \mu, \xi \rangle=\iota_{X_\xi}\omega, \quad \xi \in 
\mathfrak{t}^k$$
for the vector field $X_\xi$ on $\mathbb{C}^n$ given by the
infinitesimal action
$$X_\xi(z)=\dot{\rho}(\xi)(z), \quad z \in \mathbb{C}^n.$$
We will assume throughout
this paper hypothesis (H), i.e. the moment map is proper and
$\tau \in \mathfrak{t}^k$ is a regular value of the moment map. 
It follows from (H) that the Marsden-Weinstein quotient
$$\mathbb{C}^n//T^k=\mu^{-1}(\tau)/T^k$$
is a compact symplectic orbifold of dimension
$$\mathrm{dim}(\mathbb{C}^n//T^k)=2(n-k),$$
where the symplectic structure is
induced from the standard symplectic structure on $\mathbb{C}^n$.

\begin{rem} \label{proper}
By (\ref{moment}) the moment map is proper if and only if each 
column of the matrix $A$ has fixed sign, i.e.
for each fixed $j \in \{1,\ldots,k\}$ the sign of 
$A_{\ell j}$ is independent of $\ell \in \{1,\ldots,n\}$.
\end{rem}
We next examine which values $\tau \in \mathfrak{t}^k$ are regular. In order
to do that we first introduce some notation. 
Let $\mathscr{I}=\mathscr{I}_n$ be the set of all strictly monotone functions
$\phi \colon \{1, \ldots,n_1\} \to \{1, \ldots, n\}$ for a positive
integer $n_1 \leq n$. For each $\phi \in \mathscr{I}$ we define the
$n_1 \times k$-matrix $A_\phi$ by
$$(A_\phi)_{j \ell}:=A_{j \phi(\ell)}, \quad 1 \leq \ell \leq n_1,\,\,
1 \leq j \leq k.$$
We set
\begin{equation}\label{matrix}
\mathscr{A}:=\{A_\phi:\phi \in \mathscr {I},\,\,\mathrm{rk}(A_\phi)
<k\}.
\end{equation}
\begin{lemma}\label{irreg}
 The element of the Lie algebra
$\tau \in \mathfrak{t}^k$ is an irregular value of the moment map $\mu$, iff 
$\tau=0$ or
$\tau = A_\phi w$ where $A_\phi \in \mathscr{A}$ and the real vector
$w$ has the property that each entry is nonnegative, i.e. 
$w_j \geq 0$ for $1 \leq j \leq n_1$.
\end{lemma}
\textbf{Proof:} Assume that $\tau =A_\phi w$ where $A_\phi$ and $w$
satisfy the assumptions of the lemma. Choose a complex $n$-vector $z$
such that 
$$|z_{\phi(\ell)}|^2=2 w_\ell,\,\,1 \leq \ell \leq n_1, 
\quad z_m=0,\,\,m \in \{1, \ldots, n\}
\setminus \mathrm{im}(\phi).$$
Then
$$\mu(z)=\tau$$
and by (\ref{moment})
$$\mathrm{rk}(d\mu(z)) \leq \mathrm{rk}(A_\phi)<k$$
which shows that $\tau$ is a nonregular value of the moment map.
\\Now assume that $\tau \neq 0$ is a nonregular value of the moment map and
choose $z \in \mathbb{C}^n$ such that $\mu(z)=\tau$ and
$\mathrm{rk}(d\mu(z))<k$. Define $\phi \in \mathscr{I}$ by the property
that
$$\phi(1):=\mathrm{min}\{j \in \{1,\ldots,n\}:z_j \neq 0\},$$
$$\phi(\ell):=\mathrm{min}\{j \in \{\phi(\ell-1)+1,\ldots n\}:z_j \neq
0\},\,\,\,\,\,
2 \leq \ell \leq \#\{j \in \{1, \dots n\}:z_j \neq 0\}.$$
Define the nonnegative vector $w$ by
$$w_\ell:=\frac{1}{2}|z_{\phi(\ell)}|^2, \quad 1 \leq \ell \leq
\#\{j \in \{1, \dots n\}:z_j \neq 0\}.$$
It follows that $\tau=A_\phi w$ which proves the lemma. \hfill $\square$
\\ \\
The following Corollary of Lemma~\ref{irreg} will be used later on
to prove regularity of Givental's spaces.

\begin{cor}\label{regular}
Assume that the torus $T^k$ acts on $\mathbb{C}^n$
and $\mathbb{C}^m$ with $n\times k$-matrix $A$ respectively
$m \times k$-matrix $B$ such that for each column vector $A_\ell$ 
with $\ell \in \{1, \ldots, n\}$ of the matrix $A$ there exists a column
vector $B_{\ell'}$ of the matrix $B$ and a positive number 
$\lambda_\ell$ such that 
$$A_\ell=\lambda_\ell B_{\ell'}.$$
Suppose that $\tau \in \mathfrak{t}^k$ is an irregular value of the moment
map $\mu_A$ associated to the matrix $A$. Then $\tau$ is also an irregular
value of the moment map $\mu_B$ associated to the matrix $B$. 
\end{cor}
\textbf{Proof:} We may assume that $\tau \neq 0$. In this case
$\tau$ is an irregular nonzero value of the moment map
$\mu_A$ so that there exists by Lemma~\ref{irreg} an element $\phi \in
\mathscr{I}_n$, i.e. a strictly monotone function from
$\{1,\ldots,n_1\} \to \{1,\ldots, n\}$ for $n_1 \leq n$, and a vector
$w$ with nonnegative entries such $\mathrm{rk}(A_\phi)<k$ and
$A_\phi w=\tau$. Choose functions $\lambda \colon \{1,\dots,n\}
\to \mathbb{R}_+=\{r \in \mathbb{R}:r > 0\}$ and 
$\rho \colon \{1,\ldots,n\} \to \{1, \ldots,m\}$ such that
\begin{equation}\label{prop}
A_\ell=\lambda(\ell)B_{\rho(\ell)}
\end{equation}
which exist by the assumption of the Corollary.
Let $m_1 \leq n_1$ be the cardinality of the set
$\rho \circ \phi\{1,\ldots,n_1\}$ and define the strictly monotone
function $\tilde{\phi} \in \mathscr{I}_m$ from
$\{1, \ldots, m_1\}$ to $\{1, \ldots, m\}$ recursively by
$$\tilde{\phi}(\ell)=\min\big(\rho \circ \phi\{1,\ldots,n_1\}
\setminus \tilde{\phi}\{1,\ldots, \ell-1\}\big), \quad 
\ell \in \{1,\ldots m_1\}.$$
It follows from (\ref{prop}) that
$$\mathrm{rk}(B_{\tilde{\phi}})=\mathrm{rk}(A_\phi)<k.$$
Define for $1 \leq j \leq m_1$
$$\tilde{w}_j:=\sum_{k \in \rho^{-1}\tilde{\phi}(j)}\lambda(k)
w_{\phi^{-1}k}.$$
Then $\tilde{w}_j$ is nonnegative and we calculate
\begin{eqnarray*}
B_{\tilde{\phi}}\tilde{w}&=&\sum_{j=1}^{m_1}B_{\tilde{\phi}(j)}\tilde{w}_j\\
&=&\sum_{j=1}^{m_1}B_{\tilde{\phi}(j)}\Bigg(
\sum_{k \in \rho^{-1}\tilde{\phi}(j)}\lambda(k)w_{\phi^{-1}(k)}\Bigg)\\
&=&\sum_{j=1}^{m_1}\Bigg(\sum_{k \in \rho^{-1}\tilde{\phi}(j)}
\lambda(k)w_{\phi^{-1}(k)}B_{\rho(k)}\Bigg)\\
&=&\sum_{j=1}^{m_1}\Bigg(\sum_{k \in \rho^{-1}\tilde{\phi}(j)}
A_k w_{\phi^{-1}(k)}\Bigg)\\
&=&\sum_{i=1}^{n_1}A_{\phi(i)}w_i\\
&=&A_\phi w\\
&=&\tau.
\end{eqnarray*}
Now Lemma~\ref{irreg} implies that $\tau$ is an irregular value of
$\mu_B$. This proves the Corollary. \hfill $\square$

\section[Moment Floer homology]{Moment Floer homology}\label{moment}

Let $\mathscr{L}$ be the loop space
$$\mathscr{L}:=C^\infty(S^1,\mathbb{C}^n \times \mathfrak{t}^k),$$
where
$$S^1=\mathbb{R}/\mathbb{Z}$$
is the circle.
The gauge group
$$\mathcal{H}=C^\infty(S^1,T^k)$$
acts on $\mathscr{L}$ by
$$h_*(z,\eta)=(\rho(h)z,\eta-h^{-1}\partial_t h), \quad
h \in \mathcal{H},\,\,(z,\eta) \in \mathscr{L}.$$
The action functional for moment Floer homology 
$\mathcal{A} \colon \mathscr{L} \to \mathbb{R}$ 
is given by
$$\mathcal{A}(z,\eta):=\int_0^1 \lambda(z)(\partial_t z)+
\int_0^1\langle \mu(z(t))-\tau,\eta(t)
\rangle dt.$$
Here $\lambda$ denotes the Liouville 1-form 
$$\lambda=\sum_{i=1}^n y_i dx_i, \quad d\lambda=-\omega.$$
The first integral on the right hand side is Floer's action functional
on the space of loops in $\mathbb{C}^n$. One may think of $\eta$ in the
second integral as 
Lagrange multiplier to the constraint $\mu(z)=\tau$.

Borrowing notation from \cite{manolescu} we introduce the subgroup
$\mathcal{H}_0 \subset \mathcal{H}$ of "normalized gauge
transformations" consisting of $h \in \mathcal{H}$ for which there
exists
$\xi \in C^\infty(S^1,\mathfrak{t}^k)$ with $\int_0^1\xi dt=0$ such that
$$h(t)=\exp(\xi(t)), \quad t \in S^1.$$
The gauge group $\mathcal{H}$ can be written as
$$\mathcal{H}=\mathcal{H}_0 \times H$$
where the finite dimensional group
$$H=\{h \in \mathcal{H}: \partial_t (h^{-1}\partial_t h)=0\}$$
consists of elements $h=h_0 \cdot e^{2\pi i v t}$ with $h_0 \in T^k$
and $v \in \mathbb{Z}^k$. There is a natural isomorphism
$$H \cong T^k \times \mathbb{Z}^k$$
given by
$$h \mapsto \bigg(h(0),\frac{1}{2 \pi i}(h^{-1}\partial_t h)(0)\bigg),\quad
h \in H.$$
The action functional is invariant under $\mathcal{H}_0$
and $\mathcal{H}_0$ acts freely on $\mathscr{L}$. Moreover, 
using a gauge transformation
$h \in \mathcal{H}_0$ we can assume that $h_*\eta$ is constant. This
gauge fixing is the analogon of the Coulomb gauge fixing considered in
\cite{kronheimer-manolescu,manolescu}. In particular, the Coulomb gauge
gives a natural bijection
$$\mathscr{L}/\mathcal{H}_0 \cong C^\infty(S^1,\mathbb{C}^n)\times t^k
=:\mathscr{L}_0.$$
We consider as in \cite{kronheimer-manolescu,manolescu} the restriction
of the action functional $\mathcal{A}$ to the space $\mathscr{L}_0$
$$\mathcal{A}_0(z,\eta)=\mathcal{A}|_{\mathscr{L}_0}(z,\eta)
=\int_0^1 \lambda(z)(\partial_t z)+\Bigg\langle \eta,\int_0^1\mu(z(t))dt-\tau
\Bigg\rangle.$$
Note that for $\mathcal{A}_0$ the Lagrangian multiplier $\eta$ gives only a
constraint on the mean value of $\mu(z)$.

On $\mathscr{L}$ we have a natural $\mathcal{H}_0$-invariant
$L^2$-metric $g$ defined by integration
$$g((\hat{z}_1,\hat{\eta}_1),(\hat{z}_2,\hat{\eta}_2))
=\int_0^1 \langle\hat{z}_1,\hat{z}_2\rangle dt+\int_0^1\langle \hat{\eta}_1,
\hat{\eta}_2\rangle dt$$
where
$$(\hat{z}_1,\hat{\eta}_1),(\hat{z}_2,\hat{\eta}_2)\in 
T_{(z,\eta)}\mathscr{L}\cong \mathscr{L}, \quad (z,\eta)\in
\mathscr{L}.$$
Using the $L^2$-metric on $\mathscr{L}$ we can define two different metrics on
$\mathscr{L}_0$. The first one $g_0$ is given by the pullback of the
 inclusions $\iota \colon\mathscr{L}_0 \to\mathscr{L}$, i.e.
$g_0:=\iota^*g.$
The second one $g_1$ is the induced metric on
the quotient $\mathscr{L}_0 \cong \mathscr{L}/\mathcal{H}_0$. If
$(z,\eta)\in \mathscr{L}_0$ and 
$(\hat{z}_1,\hat{\eta}_1),(\hat{z}_2,\hat{\eta}_2) \in 
T_{(z,\eta)}\mathscr{L}_0 \cong \mathscr{L}_0$ then the first metric 
is again given by an inner product on the vector space $\mathscr{L}_0$
$$g_0((\hat{z}_1,\hat{\eta}_1),(\hat{z}_2,\hat{\eta}_2))
=\int_0^1 \langle \hat{z}_1,\hat{z}_2\rangle dt+\langle \hat{\eta}_1,
\hat{\eta}_2\rangle.$$
To give a formula for the second metric we first introduce some
notation. We define for $z \in \mathbb{C}^n$
the linear map $L_z \colon \mathfrak{t}^k \to 
T_z\mathbb{C}^n \cong \mathbb{C}^n$ by
$$L_z \xi:=X_\xi(z)=\dot{\rho}(\xi)z, \quad \xi \in \mathfrak{t}^k$$
and denote by $L^*_z \colon T_z \mathbb{C}^n \to \mathfrak{t}^k$ 
its adjoint with respect to the standard inner products
on the Lie algebra $\mathfrak{t}^k=i\mathbb{R}^n$ 
and $\mathbb{C}^n$. If 
$\xi \in \mathrm{Lie}(\mathcal{H}_0)$, i.e. 
$\xi \in C^\infty(S^1, \mathfrak{t}^k)$ and $\int_0^1 \xi dt=0$, then the 
infinitesimal action of $\xi$ on $(z,\eta) \in \mathscr{L}$ is given by
\begin{equation}\label{infgauge}
L_{(z,\eta)}\xi=(L_z \xi, -\partial_t \xi).
\end{equation}
In particular, $(\hat{z},\hat{\eta}) \in T_z \mathscr{L} \cong
\mathscr{L}$ lies orthogonal to the infinitesimal gauge action iff
$$\frac{d}{dt}\bigg(L_z^*\hat{z}+\partial_t \hat{\eta}\bigg)=0.$$
The formula for the second metric on $\mathscr{L}_0$ can now be written as
\begin{eqnarray}\label{metric2}
g_1((\hat{z}_1,\hat{\eta}_1),(\hat{z}_2,\hat{\eta}_2))&=&
g((\hat{z}_1-L_z\xi_1,\hat{\eta}_1+\partial_t \xi_1),
(\hat{z}_2-L_z\xi_2,\hat{\eta}_2+\partial_t \xi_2))\\ \nonumber
&=&\int_0^1 \big(\langle \hat{z}_1-L_z \xi_1,\hat{z}_2-L_z \xi_2\rangle 
+\langle \partial_t \xi_1,\partial_t \xi_2\rangle\big) dt\\ \nonumber
& &+\langle \hat{\eta}_1,\hat{\eta}_2 \rangle
\end{eqnarray}
where the Lie algebra elements $\xi_1,\xi_2 \in
\mathrm{Lie}(\mathcal{H}_0)$ are determined by 
$$\frac{d}{dt}\bigg(L^*_z \hat{z}_i-L_z^*L_z \xi_i+\partial_t^2 \xi_i
\bigg)=0, \quad \int_0^1\xi_idt=0, \quad i \in \{1,2\}.$$
We denote by $|| \cdot||_i$ the length of a tangent vector 
with respect to the metric $g_i$ for 
$i\in \{0,1\}$. We clearly have
$$|| \cdot||_1 \leq || \cdot ||_0.$$ 
In the following proposition we will estimate the metric $g_0$ by the
metric $g_1$. We will use this estimate later to prove that the finite
dimensional approximations of $\mathscr{L}_0$ are totally geodesic with
respect to $g_1$.
\begin{prop}\label{metest}
There exists a constant $c>0$ such that 
$$|| \cdot ||_0 \leq c\big(1+||z||_{L_2}\big)\cdot ||\cdot ||_1.$$
\end{prop}
\textbf{Proof:} Let $(z,\eta) \in \mathscr{L}_0$, 
$(\hat{z},\hat{\eta}) \in T_{(z,\eta)}\mathscr{L}_0 \cong \mathscr{L}_0$
and define $\xi \in \mathrm{Lie}(\mathcal{H}_0)$ by 
\begin{equation}\label{m1}
\frac{d}{dt}\bigg(L^*_z \hat{z}-L_z^*L_z \xi+\partial_t^2 \xi
\bigg)=0, \quad \int_0^1\xi dt=0.
\end{equation}
We calculate
\begin{eqnarray*}
0&=&\int_0^1 \langle L_z^*\hat{z}-L_z^*L_z \xi+\partial_t^2 \xi,\xi \rangle dt
\\
&=&\int_0^1\langle \hat{z},L_z\xi \rangle dt-
\int_0^1 \langle L_z \xi,L_z \xi \rangle-\int_0^1 \langle \partial_t \xi,
\partial_t \xi \rangle dt
\end{eqnarray*}
from which we obtain by the Cauchy-Schwarz inequality
$$||L_z\xi||^2_{L_2}+||\partial_t \xi||^2_{L_2}
\leq ||\hat{z}||_{L_2}||L_z \xi||_2.$$
This implies
\begin{equation}\label{m2}
||L_z\xi||_{L_2} \leq 
\frac{||L_z \xi||^2_{L_2}}{||L_z \xi||^2_{L_2}+||\partial_t
\xi||^2_{L_2}}
||\hat{z}||_{L_2}.
\end{equation}
Using the second equation in (\ref{m1}), partial integration and
H\"olders inequality we obtain as a special case of Poincar\'e's inequality 
\begin{equation}\label{m3}
||\xi||_{L_2} \leq ||\partial_t \xi||_{L_2}.
\end{equation}
The formula $L_z \xi=\dot{\rho}(\xi)z$ together with H\"olders
inequality implies that there exists a
constant $c_0>0$ such that
\begin{equation}\label{m4}
||L_z \xi||_{L_2}\leq c_0||z||_{L_2}||\xi||_{L_2}.
\end{equation}
Using inequalities (\ref{m2}), (\ref{m3}), and (\ref{m4}) we obtain the 
inequality
\begin{equation}\label{m5}
||L_z\xi||_{L_2}\leq \frac{c_0^2||z||_{L_2}^2}{c_0^2||z||_{L_2}^2+1}||
\hat{z}||_{L_2}.
\end{equation}
We now estimate using (\ref{metric2}) and (\ref{m5})
\begin{eqnarray*}
||(\hat{z},\hat{\eta})||_1^2 &\geq&
||\hat{z}-L_z\xi||^2_{L_2}+||\hat{\eta}||^2_{L_2}\\
&\geq&||\hat{z}||^2_{L_2}-||L_z\xi||^2_{L_2}+||\hat{\eta}||^2_{L_2}\\
&\geq&||\hat{z}||^2_{L_2}-\Bigg
(\frac{c_0^2||z||_{L_2}}{c_0^2||z||^2_{L_2}+1}\Bigg)^2
||\hat{z}||^2_{L_2}+||\hat{\eta}||^2_{L_2}\\
&\geq&||\hat{z}||^2_{L_2}-
\frac{c_0^2||z||_{L_2}}{c_0^2||z||^2_{L_2}+1}
||\hat{z}||^2_{L_2}+||\hat{\eta}||^2_{L_2}\\
&\geq&\frac{1}{c_0^2||z||^2_{L_2}+1}||(\hat{z},\hat{\eta})||_0.
\end{eqnarray*}
The Proposition follows now with $c=\max\{c_0,1\}$. \hfill $\square$
\\ \\
We denote by $\nabla_0\mathcal{A}_0$ the gradient of $\mathcal{A}_0$
on $\mathscr{L}_0$ with respect to the metric $g_0$ and by
$\nabla_1 \mathcal{A}_0$ the gradient with respect to the metric $g_1$.
We will use the notation
$$\bar{\mu}(z):=\int_0^1 \mu(z(t))dt.$$

\begin{prop} \label{gradienten}
The two gradients are given by
\begin{equation}\label{grad1}
\nabla_0\mathcal{A}_0(z,\eta)=
(i\partial_t z+iL_z \eta,\bar{\mu}(z)-\tau)
\end{equation}
and
\begin{equation}\label{grad2}
\nabla_1 \mathcal{A}_0(z,\eta)=
(i\partial_t z+iL_z \eta+L_z \xi, \bar{\mu}(z) -\tau)
\end{equation}
where $\xi \in \mathrm{Lie}(\mathcal{H}_0)$ is determined by
\begin{equation}\label{gradxi}
\partial_t \xi(t)=\mu(z(t))-\bar{\mu}(z), \quad \int_0^1\xi dt=0.
\end{equation}
\end{prop}
\textbf{Proof: } To determine the first gradient we compute 
for $(\hat{z},\hat{\eta})\in T_{(z,\eta)}\mathscr{L}_0\cong \mathscr{L}_0$
\begin{eqnarray*}
d\mathcal{A}_0(z,\eta)(\hat{z},\hat{\eta})
&=&\int_0^1 d\lambda(z)(\hat{z},\partial_t z)
+\langle \eta, \int_0^1 d\mu(z(t))\hat{z}(t)dt \rangle
+\langle \hat{\eta},\bar{\mu}(z)-\tau \rangle\\
&=&-\int_0^1\omega(\hat{z},\partial_t z)+\int_0^1 \omega(X_\eta(z),\hat{z})
+\langle \hat{\eta},\bar{\mu}(z)-\tau \rangle\\
&=&\langle (i\partial_t z+iL_z\eta, \bar{\mu}(z)-\tau),(\hat{z},\hat{\eta})
\rangle_1
\end{eqnarray*}
which implies (\ref{grad1}).

To compute the second gradient we make use of the fact that the
action functional $\mathcal{A}$ on $\mathscr{L}$ is invariant under
$\mathcal{H}_0$ and hence is gradient $\nabla\mathcal{A}$
is orthogonal to the infinitesimal
gauge action. For 
$\zeta=(\hat{z},\hat{\eta}) \in T_{(z,\eta)}\mathscr{L}_0$ we
denote by $\xi_\zeta \in \mathrm{Lie}(\mathcal{H}_0)$ the Lie algebra
element which is determined by
$$\frac{d}{dt}\Bigg(L_z^*\hat{z}-L_z^*L_z\xi_\zeta
+\partial^2_t\xi_\zeta\Bigg)=0,\quad \int_0^1\xi_\zeta dt=0.$$
We will abbreviate
$$\xi=\xi_{\nabla_1\mathcal{A}_0(z,\eta)}.$$
We compute for the metric $\langle \cdot,\cdot \rangle$ on
$\mathcal{L}$
\begin{eqnarray*}
\langle\nabla_1\mathcal{A}_0(z,\eta)-
L_{(z,\eta)}\xi,
\zeta-L_{(z,\eta)}\xi_{\zeta}\rangle
&=&\langle \nabla_2\mathcal{A}_0(z,\eta),\zeta\rangle_2\\
&=&d\mathcal{A}_0(z,\eta) \zeta\\
&=&d\mathcal{A}(z,\eta)\zeta\\
&=&\langle \nabla\mathcal{A}(z,\eta),\zeta\rangle\\
&=&\langle \nabla \mathcal{A}(z,\eta),\zeta-L_{(z,\eta)}\xi_\zeta\rangle
\end{eqnarray*}
which implies
\begin{equation}\label{gr20}
\nabla_1\mathcal{A}_0(z,\eta)-
L_{(z,\eta)}\xi=
\nabla \mathcal{A}(z,\eta).
\end{equation}
By a computation similar to the computation of $\nabla_0\mathcal{A}_0$
one gets
\begin{equation}\label{grad0}
\nabla \mathcal{A}(z,\eta)=(i\partial_t z+iL_z\eta,\mu(z)-\tau).
\end{equation}
Using (\ref{infgauge}),(\ref{gr20}), and (\ref{grad0}) we conclude that
$\mu(z)-\tau-\partial_t \xi$ is independent of the
$t$-variable. Observing that $\xi(0)=\xi(1)$ we obtain the formula
(\ref{gradxi}) for $\xi$. Using again (\ref{gr20}) and (\ref{grad0})
we are able to deduce the formula (\ref{grad2}) for the second gradient.
This completes the proof of the proposition.\hfill $\square$

\section[Finite dimensional approximation]{Finite dimensional
 approximation}\label{finite}

Formulas (\ref{grad1}) and (\ref{grad2}) show that the gradients of
Floer's action functional are a zero'th order perturbation of the 
Cauchy-Riemann operator. Following 
\cite{furuta, kronheimer-manolescu, manolescu} we will approximate the
loop space by the eigenspaces of the first order part of the gradient, which 
leads in our case to Fourierapproximation.

For $1 \leq j \leq n$ let $m^-_j \leq m^+_j$ be integers. We consider 
the finite dimensional complex vector space
$$V=V_{\{m^-_j,m^+_j\}_{1 \leq j \leq n}}$$
consisting of finite Fourier series
$z=(z_1, \ldots, z_n) \in C^\infty(S^1,\mathbb{C}^n)$ 
for which there exist complex numbers $z_{jm}$ where $1 \leq j \leq n$
and $m^-_j \leq m \leq m^+_j$ such that
$$z_j(t)=\sum_{m=m^-_j}^{m^+_j}z_{jm}e^{2\pi i m t}, \quad j \in 
\{1,\ldots,n\},\,\,t \in S^1.$$
The action of the torus $T^k$ on $\mathbb{C}^n$ induces an action of
$T^k$ on $V$ by pointwise multiplication
$$(\gamma z)(t):=\gamma(z(t)), \quad \gamma \in T^k,\,\,z \in V,\,\,t \in S^1.$$
Setting
$$N:=\sum_{j=1}^n (m^+_j-m^-_j+1)$$
we will identify in the following $V$ with $\mathbb{C}^N$ using the 
map
$$z \mapsto (z_{1m_1^-},\ldots ,z_{1m_1^+},z_{2m_2^-},\ldots
,z_{nm^+_n}).$$
Denote by $\mu_V$ the moment map associated by (\ref{moment}) to the
action of $T^k$ on $V$. The restriction of the action functional 
$\mathcal{A}_0$ to the finite dimensional space $V \times
\mathfrak{t}^k$ is given by the formula
\begin{equation}\label{finac}
\mathcal{A}_0|_{V \times \mathfrak{t}^k}(z,\eta)
=\int_0^1 \lambda(z)(\partial_t z)+\langle \eta,\mu_V(z)-\tau \rangle,\quad
z \in V,\,\,\eta \in \mathfrak{t}^k.
\end{equation}
Denote by $A_V$ the matrix corresponding to the moment map $\mu_V$ defined in
(\ref{matrix}). To each column vector of $A_V$ corresponds a column
vector of $A$, more precisely let 
$\rho \colon \{1,\ldots,N\} \to \{1,\ldots,n\}$
be defined by
$$\rho(k):=\min\{j: \sum_{i=1}^j(m^+_i-m^-_i+1)\geq k\}, \quad 1\leq
k\leq N$$
then
$$(A_V)_k=A_{\rho(k)}, \quad 1 \leq k\leq N.$$
Using this formula together with Remark~\ref{proper} and 
Corollary~\ref{regular} we obtain the following proposition.
\begin{prop}\label{neureg}
Assume (H), i.e. the moment map 
$\mu$ is proper and $\tau$ is a regular value of $\mu$. Then $\mu_V$
is also proper and $\tau$ is also a regular value of $\mu_V$. 
\end{prop}
\begin{rem}
It follows from Proposition~\ref{neureg} that the spaces
$\mu_V^{-1}(\tau)$ are compact manifolds. Their quotients under
the torus action $\mu_V^{-1}(\tau)/T^k$ are therefore compact symplectic
orbifolds. They are the toric map spaces introduced by Givental in 
\cite{givental2}
as an approximation to the
moduli spaces of Floer homology for toric orbifolds, see also
\cite{givental, iritani, vlassopoulos}.
\end{rem}
We next examine the metrics $g_0$ and $g_1$ of $\mathscr{L}_0$
on the finite dimensional submanifold $V \times \mathfrak{t}^k$.
The metric $g_0$ equals the metric induced by the standard scalar
product on $V \times \mathfrak{t}^k \cong \mathbb{C}^N \times
\mathfrak{t}^k$. It follows from (\ref{metric2}) that $g_1$
on $V \times \mathfrak{t}^k$ is given by a product metric
$$g_1=g_V\oplus g_{\mathfrak{t}^k}$$
where $g_{\mathfrak{t}^k}$ is induced by the standard scalar product on
$\mathfrak{t}^k$ and $g_V$ is a Riemannian metric on $V$. 
\begin{prop}\label{complete}
The Riemannian metric $g_V$ on $V$ is geodesically complete.
\end{prop}
\textbf{Proof: }We have to show that the length of each path
$\gamma \colon [a,b) \to V\cong \mathbb{C}^N$ for $a,b \in \mathbb{R}$
such that $\gamma(s)$ tends to infinity as $s$ goes to $b$ has infinite
length with respect to $g_V$. Define the path
$\tilde{\gamma} \colon [0,\infty) \to \mathbb{C}^N$ as the
reparametrisation of $\gamma$, which is parametrized
by arclength with respect to the standard metric on $\mathbb{C}^N$. 
We estimate using Proposition~\ref{metest}
\begin{eqnarray*}
\mathrm{length}(\gamma)&=&
\int_a^b\sqrt{g_V(\gamma(t))(\partial_t \gamma(t),\partial_t \gamma(t))}dt\\
&=&\int_0^\infty \sqrt{g_V(\tilde{\gamma}(t))(\partial_t \tilde{\gamma}(t),
\partial_t \tilde{\gamma}(t))}dt\\
&\geq&\int_0^\infty\frac{1}{c(1+||\tilde{\gamma}(t)||)}
||\partial_t\tilde{\gamma}(t)|| dt\\
&\geq&\int_0^\infty\frac{1}{c(1+||\tilde{\gamma}(0)||+t)}dt\\
&=&\infty.
\end{eqnarray*}
This proves the proposition. \hfill $\square$

\section[Conley index]{Conley index}\label{conley}

Assume that $M$ is a finite dimensional manifold and $\phi$ is a flow on
$M$, i.e. a continuous map $\phi \colon M \times \mathbb{R} \to M,
(x,t) \mapsto \phi^t (x)$, which satisfies $\phi^0=\mathrm{id}$ and
$\phi^s \circ \phi^t=\phi^{s+t}$. For a subset $A \subset M$ we set
$$\mathrm{Inv}_\phi A=\{x \in A: \phi^t(x) \in A,\,\,\forall\,\,t \in
\mathbb{R}\}.$$
A subset $S \subset M$ is called an \textbf{isolated invariant
set} if there exists an isolating neighbourhood $A$ for $S$,
i.e. a compact set $A$ such that
$$S =\mathrm{Inv}_\phi(A) \subset \mathrm{int}(A).$$
In particular, $S$ is compact and $S=\mathrm{Inv}_\phi (S)$.
The Conley index $I(\phi,S)$ attributes to every isolated invariant set $S$
of the flow $\phi$ an isomorphism class of a pointed topological space,
i.e. a topological space with a distinguished based point. We refer the
reader to \cite{conley, salamon} for an introduction into Conley index
theory and to \cite{manolescu} for a short survey. We will use in this
paper the following two properties of the Conley index.
\\ \\
\textbf{Property 1 (Invariance)} \emph{Assume that $A$ is an isolating 
neighbourhood for $S_r=\mathrm{Inv}_{\phi_r} A$ for a continuous family of
flows $\phi_r, r\in [0,1]$, then $I(\phi_0,S_0)=I(\phi_1,S_1)$.}
\\ \\
To formulate the second property we have first to introduce some
notation. Assume that $E \to M$ is a vector bundle. Choose a bundle
metric $g$ on $E$ and denote by $B_gE$ and $S_g E$ the unit ball - 
respectively the unit sphere bundle - bundle over M. Then the Thom space
of $E$ is defined as the pointed topological space
$$T_gE:=B_gE /S_gE:=((B_gE \setminus S_gE \cup[S_gE]),[S_gE])$$
obtained by collapsing $S_g E$ to a single point denoted by $[S_g E]$.
A set $U \subset T_gE$ is open if either $U$ is open in $B_g E$
and $U \cap S_g E=\emptyset$ or the set 
$(U \cap (B_g E \setminus S_g E)) \cup S_g E$ is open in $B_g E$.
The pointed isomorphism class of $T_g E$ is independent of the choice of
the bundle metric $g$ and will be denoted by $TE$. \\
Now assume that $\phi$ is a gradient flow on the manifold $M$, i.e. 
there exists a function $f \in C^\infty(M,\mathbb{R})$ and a Riemannian
metric $g$ on $M$ such that
$$\frac{d}{dt}\phi^t(x)=-\nabla_g f(\phi^t(x)), \quad x\in M.$$
Now assume that the subset $S \subset M$ is a critical
submanifold of Morse-Bott type of the function $f$, i.e. $S$ is a 
submanifold of $M$ and for each $x \in S$
$$T_x S=\mathrm{ker}H_f(x)$$
where $H_f$ denotes the Hessian of $f$. Under this assumption there
is a natural splitting of the normal bundle $NS$ of $S$ in $M$
$$NS=N_f^+S \oplus N_f^- S$$
where for $x \in S$, $(N_f^+S)_x$ is spanned by the eigenvectors of the Hessian
$H_f(x)$ to positive eigenvalues and $(N_f^-S)_x$ is spanned by the 
eigenvectors of $H_f(x)$ to negative eigenvalues. 
\\
We are now in the position to formulate the second property about the Conley
index we will need.
\\ \\
\textbf{Property 2 (Nontriviality)} \emph{Assume that $\phi$ is a gradient flow
with respect to the function $f \in C^\infty(M,\mathbb{R})$ and
$S$ is a compact critical submanifold of Morse-Bott type of $f$.
Then 
$$I(\phi,S)=T(N_f^-S).$$}
\\
In this section we are interested in some special class of isolated
invariant sets. 
Assume that $X \in \Gamma(TM)$ is a vector field on $M$.
For a flow line of $X$, i.e. a map $y \in C^\infty(\mathbb{R},M)$
satisfying
$$\partial_s y(s)=X(y(s)), \quad s \in \mathbb{R}$$
we define its $\omega$-limit sets as
$$\omega^+(y)=\bigcap_{t=0}^\infty\mathrm{cl}
\Bigg(\bigcup_{s \geq t}y(s)\Bigg),\quad
\omega^-(y)=\bigcap_{t=0}^\infty\mathrm{cl}
\Bigg(\bigcup_{s \geq t}y(-s)\Bigg).$$
We denote by $S_X$ the subset of $M$ consisting of the traces of the
flow lines of the vector field $X$ both of whose $\omega$-limit sets
are nonempty, i.e.
$$S_X:=\{y(\sigma): \sigma \in \mathbb{R},\,\,y \in 
C^\infty(\mathbb{R},M),\,\,\partial_s y=X(y),\,\,
\omega^\pm(y) \neq \emptyset\}.$$
\begin{rem}
We point out that even for gradient vector fields the $\omega$-limit
sets may consist of more than one point, see \cite{palis-demelo}. 
However, if the critical set of a function consists of critical
manifolds of Morse-Bott type then the assumption of
nonempty $\omega$-limit sets
implies that the flow line converges at both ends. In our applications to
the finite dimensional approximations of the action functional 
$\mathcal{A}_0$ we always have Morse-Bott situations, so that for
these cases the set $S_X$ may be defined as the traces of flow lines
which converge at both ends.
\end{rem}
Assume that $S_X$ is compact. 
If the manifold $M$ is not compact then
the vector field $X$ may not generate a well-defined global flow. To
overcome this problem choose compact
subsets $A,A' \subset M$ such that
$S_X \subset \mathrm{int}(A) \subset A \subset A'$ and a smooth
cutoff function $\beta$ such that $\beta|_{A}=1$ and $\beta|_{M
\setminus A'}=0$. Denote by $\phi_X=\phi_{X,\beta}$ the flow of the 
gradient of the
vector field $\beta \cdot X$. Since $\beta \cdot X$ has compact support the
flow $\phi_X$ is well defined. Moreover, the set $S_X$ is an isolated
invariant set with isolating neighbourhood $A$. We now set
$$I(X)=I(\phi_X,S_X).$$
It remains to check that $I(X)$ is indeed well-defined, i.e. independent
of the choice of the cutoff function $\beta$. This follows readily from
Property 1 of the Conley index. Indeed, assume that $\beta_1$
and $\beta_2$ are two cutoff functions satisfying the conditions above
for compact sets $A_1,A_1'$ respectively $A_2,A_2'$ then for the family
of flows
$\{\phi_r\}, r \in [0,1]$ of the vector field 
$(r \beta_1+(1-r)\beta_2)\cdot X$ the set $A_1 \cap A_2$ is an
isolating neighbourhood of $S_X$ and hence the Conley indices for $\phi_0$
and $\phi_1$ agree by Property 1. 
\\
We are now in position to state the main theorem of this section.

\begin{thm}\label{palais}
Assume that $(M,g)$ is a geodesically complete Riemannian manifold,
$f \in C^\infty(M,\mathbb{R})$ and $h \in C^\infty(M,\mathbb{R}^k)$
for $k \in \mathbb{N}$ are smooth functions such that $h$ is proper
and $0$ is a regular value of $h$. Define 
$F \in C^\infty(M \times \mathbb{R}^k,\mathbb{R})$ by
$$F(x,\lambda)=f(x)+\langle h(x),\lambda \rangle,\quad
x \in M,\,\,\lambda \in \mathbb{R}^k$$
and denote by $\nabla F$ its gradient with respect to the product
metric $g \oplus g_{\mathbb{R}^k}$ on $M \times \mathbb{R}^k$ where 
$g_{\mathbb{R}^k}$
is the metric on $\mathbb{R}^k$ obtained by the standard scalar product.
Then $S_{-\nabla F}$ is compact and
\begin{equation}\label{suco}
I(-\nabla F)=T(N_M h^{-1}(0))
\end{equation}
where $N_Mh^{-1}(0)$ denotes the normal bundle of $h^{-1}(0)$ in $M$.
In particular, $I(-\nabla F)$ does not depend on $f$.
\end{thm}

\begin{fed}\label{tame}
A vector field $X \in \Gamma (TM)$ on a manifold $M$
is called \emph{\textbf{tame}} if
there exists a compact set $K \subset M$, a constant
$\epsilon>0$,a geodesically complete Riemannian
metric $g$ on $M$, and a smooth function
$F \in C^\infty(M,\mathbb{R})$ satisfying the following properties
\begin{description}
 \item[(i)] $X$ does not vanish outside of $K$, i.e.
  $$X(x)=0\,\, \Longrightarrow\,\, x \in K.$$
 \item[(ii)] $F$ is a Lyapunov function, i.e.
  $$dF(x)X(x) \leq 0, \quad x \in M.$$
 \item[(iii)] $X$ satisfies the following Palais-Smale condition
  $$dF(x)X(x)\leq -\epsilon \sqrt{g(X(x),X(x))}, \quad x \in M \setminus K.$$
\end{description}
A smooth family of vector fields $X_r \in \Gamma (TM)$ for $r \in [0,1]$
is called tame if there exist a compact set $K \subset M$,
a constant $\epsilon>0$, and a geodesically complete Riemannian metric
$g$ independent of $r$, and a family of smooth
functions $F_r \in C^\infty(M,\mathbb{R})$ which are uniformly bounded
on $K$ satisfying conditions
(i),(ii), and (iii) above for every $r \in [0,1]$.
\end{fed}

\begin{prop}\label{lya}
Assume that $X \in \Gamma (TM)$ is a tame vector field. 
Then $S_X$ is compact. Moreover,
if $X_r \in \Gamma(TM)$ for $r \in [0,1]$ is a tame family of vector
fields then $I(X_r)$ does not depend on $r \in [0,1]$.
\end{prop}
\textbf{Proof:} Assume that $y \in C^\infty(\mathbb{R},M)$ is a flow line of $X$, i.e.
$$\partial_s y=X(y),$$
such that $\omega^\pm(y)\neq \emptyset$. We show that there exists a
compact set $K' \subset M$ independent of $y$ such that the trace of
$y$ is entirely contained in $K'$, i.e.
$$y(\sigma) \in K', \quad \sigma \in \mathbb{R}.$$
Choose $y^\pm \in \omega^\pm(y)$. It follows from the Lyapunov property 
(ii) in Definition~\ref{tame} that $dF(y^\pm)X(y^\pm)=0$. Now properties
(i) and (iii) imply that $y^\pm \in K$. Since $K$ is
compact there exists a constant $m>0$ independent of $y$ such that
$$|F(y^\pm)| \leq m.$$
Using again (ii) we conclude that
$$|F(y(\sigma))| \leq m, \quad \sigma \in \mathbb{R}.$$
For $s \in \mathbb{R}$ define
$$\tau(s):=\inf_{s' \geq s}\{y(s') \in K\}.$$
Here we use the convention that the infimum of the empty set equals
infinity.
We now estimate the distance from $y(s)$ to $K$ by
\begin{eqnarray*}
d(y(s),K) &\leq& \int_s^{\tau(s)} ||\dot{y}(\sigma)|| d\sigma\\
&=&\int_s^{\tau(s)}||X(y(\sigma))||d \sigma\\
&\leq& -\frac{1}{\epsilon}
\int_s^{\tau(s)}dF(y(\sigma))X(y(\sigma))d\sigma\\
&\leq& -\frac{1}{\epsilon}
\int_{-\infty}^\infty dF(y(\sigma))X(y(\sigma))d\sigma\\
&=&-\frac{1}{\epsilon}\int_{-\infty}^\infty 
\frac{d}{d\sigma}F(y(\sigma))d\sigma\\
&=&\lim_{s \to \infty}\frac{1}{\epsilon}\Bigg(F(y(-s))-F(y(s))\Bigg)\\
&\leq&\frac{2m}{\epsilon}.
\end{eqnarray*}
We now define
$$K'=\{x \in M:d(x,K)\leq 2m /\epsilon\}.$$
Since the metric $g$ is geodesically complete the set $K'$ is
compact. Hence the set $S_X$ is compact as a closed subset of the
compact set $K'$.
\\
If $X_r \in \Gamma(TM)$ for $r \in [0,1]$ is a tame family of vector
fields, we define $K'$ as above, by choosing as $m$ the uniform bound
of the family of functions $F_r$ on $K$. The above estimate then shows
that $S_{X_r}$ is contained in $K$ for every $r \in [0,1]$.
Now choose compact sets $A,A'$ such that $K \subset \mathrm{int}(A)
\subset A \subset A'$ and a smooth cutoff function $\beta$ such 
that $\beta|_{A}=1$ and $\beta_{M \setminus A'}=0$. Define
$\phi_r$ to be the flow of the compactly supported vector field
$\beta \cdot X_r$ for $r \in [0,1]$. Then $A$ is an isolating
neighbourhood for $S_{X_r}=\mathrm{Inv}_{\phi_r} A$ and the
invariance of the Conley indices $I(X_r)$ follows from Property 1. 
This proves the Proposition. \hfill $\square$ 

\begin{lemma}\label{compact}
Assume that $(M,g)$ is a geodesically complete
Riemannian manifold and let
$f \in C^\infty(M,\mathbb{R})$ and $h \in C^\infty(M,\mathbb{R}^k)$
for $k \in \mathbb{N}$ satisfy the assumptions of Theorem~\ref{palais}.
For $r \in [0,1]$
define the functions $F_r \in C^\infty(M \times \mathbb{R}^k)$ by
$$F_r(x,\lambda):=rf(x)+\langle \lambda,h(x)\rangle,\quad x \in M,\,\lambda \in
\mathbb{R}^k.$$
Denote by $\nabla F_r$ the gradient of $F_r$ with respect to the product
metric $g \oplus g_{\mathbb{R}^k}$. Then the family of vector fields
$X_r=-\nabla F_r$ is tame.
\end{lemma}
\textbf{Proof: }Since $h$ is proper and $0$ is a regular value of $h$
there exists $\epsilon>0$ such that every 
$\lambda \in B_\epsilon=\{\mu \in \mathbb{R}^k: ||\mu|| \leq \epsilon\}$
is a regular value of $h$. Consider the continuous function
$$\rho \colon h^{-1}(B_\epsilon) \to \mathbb{R}, \quad
x \mapsto \mathrm{min}\bigg\{||\sum_{i=1}^k \lambda_i \nabla h_i(x)||:
\lambda \in \mathbb{R}^k,\,\,||\lambda||=1\bigg\}$$
where $h_i$ for $1\leq i\leq k$ denotes the $i$-th component of the
function $h \in C^\infty(M,\mathbb{R}^k)$. Since $dh(x)$ is surjective
for every $x \in h^{-1}(B_\epsilon)$
it follows that $\rho(x)>0$ and hence
$$\delta:=\mathrm{min}\{\rho(x): x\in h^{-1}(B_\epsilon)\}>0.$$
Set
$$c:=\mathrm{max}\{||\nabla f(x)||: x \in h^{-1}(B_\epsilon)\}$$
and define
$$K:=\bigg\{(x,\lambda) \in M \times \mathbb{R}^k:
x \in h^{-1}(B_\epsilon),\,\,||\lambda||\leq \frac{c+\epsilon}{\delta}
\bigg\}.$$ 
The set $K$ is compact and we claim
\begin{equation}\label{smale}
||\nabla F_r (x,\lambda)|| >\epsilon, \quad r \in [0,1],\,\,
(x,\lambda) \in M \times \mathbb{R}^k \setminus K.
\end{equation}
In order to prove (\ref{smale}) we first compute the gradient of $F_r$
$$\nabla F_r(x,\lambda)=\left(\begin{array}{c}
r \cdot \nabla f(x)+\sum_{i=1}^k\lambda_i \cdot \nabla h_i(x)\\
h(x)
\end{array}\right).$$
We first consider the case where $x \notin h^{-1}(B_\epsilon)$. 
Then 
$$||\nabla F_r(x,\lambda)|| \geq ||h(x)||>\epsilon$$
and (\ref{smale}) holds. Now consider the case where 
$x \in h^{-1}(B_\epsilon)$ but 
$||\lambda|| >(c+\epsilon)/\delta$. We estimate for every $r \in [0,1]$
\begin{eqnarray*}
||\nabla F_r(x,\lambda)|| &\geq&
||r \cdot \nabla f(x)+\sum_{i=1}^k \lambda_i \cdot \nabla h_i(x)||\\
&>&||\sum_{i=1}^k \lambda_i \cdot \nabla h_i(x)||
-||r \cdot \nabla f(x)||\\
&\geq&\delta\cdot\bigg(\frac{c+\epsilon}{\delta}\bigg)-c\\
&=&\epsilon
\end{eqnarray*}
which proves (\ref{smale}). The lemma now follows for the obvious
choices of the ingredients in the definition of tame family, namely
choose $K$ as the compact set in $M \times \mathbb{R}$, choose
$\epsilon$ as the positive constant, choose $g \oplus g_0$ as the 
geodesically complete Riemannian metric, and choose $F_r$ itself as the
family of functions. \hfill $\square$
\\ \\
\textbf{Proof of Theorem~\ref{palais}: }
We consider the homotopy of functions 
$$F_r(x,\lambda):=rf(x)+\langle \lambda, h(x)\rangle,\quad
x \in M,\,\,\lambda \in \mathbb{R}^k$$
for $r \in [0,1]$. Note that $F_1=F$. It follows from
Proposition~\ref{lya} and Lemma~\ref{compact} that
\begin{equation}\label{c1}
I(-\nabla F)=I(-\nabla F_0).
\end{equation}
It remains to compute $I(-\nabla F_0)$.
Using the formula for the gradient of 
$F_0(x,\lambda)=\langle \lambda, h(x)\rangle$
$$\nabla F_0(x,\lambda)=\left(\begin{array}{c}
\sum_{i=1}^k\lambda_i \cdot \nabla h_i(x)\\
h(x)
\end{array}\right)$$ 
and the fact that $0$ is a regular value of $h$ we conclude that the
critical set $C$ of $F_0$ is a manifold and consists of
$$C=\{(x,0):h(x)=0\}.$$
In particular, $F_0$ is constant equal $0$ on $C$ and hence all flow lines
whose limits exist on both sides are constant paths in $C$. The Hessian
$H_{F_0}$ at a point $(x,0) \in C$ is given by
\begin{equation}\label{Hess}
H_{F_0}(x,0)=\left(\begin{array}{cc}
0 & dh(x)^*\\
dh(x) & 0
\end{array}\right)
\end{equation}
where the adjoint $dh(x)^*$ is taken with respect to the metric $g_x$
on $T_x M$ and the standard scalar product on $T_{h(x)}\mathbb{R}^k
\cong \mathbb{R}^k$.
Using the assumption that $0$ is a regular value of $h$ we conclude that
$$\mathrm{rk}(H_{F_0}(x,0))=2k=\mathrm{codim}_{M \times \mathrm{R}^k}C.$$
This shows that $C$ is a critical submanifold of Morse-Bott type. 
The normal bundle $NC$ of $C$ in $M \times \mathbb{R}^k$ splits
into the sum 
$$NC=N^+_{F_0}C
\oplus N^-_{F_0}C$$ 
given by the positive, respectively
negative, eigenvalues of the Hessian $H_{F_0}$ and we conclude from
property 2 of the Conley index that
\begin{equation}\label{c2}
I(-\nabla F_0)=T(N^-_{F_0}C).
\end{equation}
In view of Lemma~\ref{bundle} below we conclude from (\ref{c1}) and
(\ref{c2}) that (\ref{suco}) holds, which proves the theorem.
\hfill $\square$

\begin{lemma}\label{bundle}
The normal bundle $N_M h^{-1}(0)$ of $h^{-1}(0)$ in $M$ and the bundle
$N^-_{F_0}C$ are isomorphic vector bundles.
\end{lemma}
There is a natural diffeomorphism from $C$ to
$h^{-1}(0)$ given by the map $(x,0) \mapsto x$. Using (\ref{Hess})
we conclude that $N^-_{F_0}C$ at $(x,0)$
is spanned by the vectors
$$w_v:=\left(\begin{array}{c}
v \\
-\frac{1}{\lambda}dh(x)v
\end{array}\right)$$
where $v$ is an eigenvector of $dh(x)^* dh(x)$ to the eigenvalue
$\lambda^2$, i.e.
$$dh(x)^* dh(x)v=\lambda^2v.$$
We can now identify the vector bundle $N^-_{F_0}C$
over $(x,0)$ with the normal bundle $N_M h^{-1}(0)$ over $x$
by the linear extension of the map
$$w_v \mapsto v.$$
This finishes the proof of the lemma. \hfill $\square$
\\ \\
For a finite dimensional approximation $V$ of the loop space
$C^\infty(S^1,\mathbb{C}^n)$ we define
$I^{g_i}_V$ for $i \in \{0,1\}$ as the Conley index for the gradient of 
$\mathcal{A}_0|_{V \times \mathfrak{t}^k}$ with respect to the metric
$g_i$. The proof of Theorem A is now straightforward.
\\ \\ 
\textbf{Proof of Theorem A:} It follows from formula (\ref{finac})
that the action functional $\mathcal{A}_0|_{V \times \mathfrak{t}^k}$
is the sum of the restriction of Floer's action functional 
to $V$ and a Lagrange multiplier to the condition
$\mu_V^{-1}(\tau)$. It follows from Proposition~\ref{neureg} that
the map $\mu_V-\tau$ is proper and $0$ is a regular value of it.
The metric $g_0$ is obviously geodesically complete on $V$ and for
the metric $g_1$ this follows from Proposition~\ref{complete}. Theorem A
follows now from Theorem~\ref{palais}. \hfill $\square$

\section[Cobordism]{Cobordism}\label{cobordism}

In this section we prove Theorem B and Theorem C. Theorem B follows from
the following Theorem.

\begin{thm} \label{embedding}
Let 
$(z,\eta) \in C^\infty(\mathbb{R}\times S^1,\mathbb{C}^n) \times
C^\infty(\mathbb{R},\mathfrak{t}^k )$ 
be a gradient flow line of $\mathcal{A}_0$
with respect to the metric $g_0$ on $\mathscr{L}_0$, i.e.
$$\partial_s (z,\eta)(s,\cdot)=-\nabla_0 \mathcal{A}_0((z,\eta)(s,\cdot)),
\quad s \in \mathbb{R}.$$
Assume the flow line converges on both ends, i.e. 
$$\lim_{s \to \pm \infty}(z,\eta)(s)=(z^\pm,\eta^\pm) \in \mathscr{L}_0,$$
where the limits are taken with respect to the $C^\infty$-topology,
exist. Denote for $1 \leq j \leq n$
$$m_j^-:=\min_{n \in \mathbb{Z}}\{2 \pi n \geq \sum_{r=1}^k A_{jr}\eta^-_r\}, 
\quad
m_j^+:=\max_{n \in \mathbb{Z}}\{2 \pi n \leq \sum_{r=1}^k A_{jr}\eta^+_r
\},$$
where $A$ is as usual the $n \times k$-matrix which defines the
action of the torus on $\mathbb{C}^n$.
Then the trace of the flow line $(z,\eta)$ is contained in the finite
dimensional subspace $V=V_{\{m_j^-,m_j^+\}_{1 \leq j \leq n}}$ of
$\mathscr{L}_0$, i.e.
$$(z,\eta)(s) \in V , \quad \forall \,\, s \in \mathbb{R}.$$
\end{thm}
\textbf{Proof:} The proof is based on Fourierapproximation. Let
$$z_j(s,t)=\sum_{m=-\infty}^\infty z_{jm}(s)e^{2 \pi i m t}, \quad
1 \leq j \leq n.$$
The formula (\ref{grad1}) shows
that for $j \in \{1,\dots, n\}$ and $m \in \mathbb{Z}$
$$\partial_s z_{jm}(s)+\Bigg(\sum_{r=1}^k A_{jr}\eta_r(s)-2\pi m
\Bigg)z_{jm}(s)=0.$$
Using $\lim_{s \to \pm \infty} \partial_s z_{jm}(s)=0$ we conclude that
$z_{jm}$ vanishes identically if $m$ is not contained in
$\{m_j^-, \ldots, m_j^+\}$. \hfill $\square$
\\ \\
Our next aim is to prove Theorem C. To motivate our construction 
of the homotopy of metrics $g_r$ on
$\mathscr{L}_0$ for $r \in [0,1]$
it is useful to  consider the family of action functionals
$$\mathcal{A}^r \colon \mathscr{L}=
C^\infty(S^1,\mathbb{C}^n \times \mathfrak{t}^k) \to \mathbb{R},
\quad r \in [0,1]$$
given by
$$\mathcal{A}^r (z,\eta):=\int_0^1 \lambda(z)(\partial_t z)
+\int_0^1 \langle (1-r)\mu(z(t))+r \bar{\mu}(z)-\tau,\eta(t)\rangle
dt.$$
Note that
$$\mathcal{A}^0=\mathcal{A}, \quad \mathcal{A}^r|_{\mathscr{L}_0}
=\mathcal{A}_0,\,\, r\in [0,1].$$
The action functionals $\mathcal{A}^r$ for $r \in [0,1]$
are invariant under the usual
action of the finite dimensional group $H$
and for $r \in (0,1]$ they are also invariant
under the following deformed action of the gauge group
$\mathcal{H}_0$ on $\mathscr{L}$
$$h_{*_r}(z,\eta)=(\rho(h)z,\eta-\frac{1}{r}h^{-1}\partial_t h),
\quad h \in \mathcal{H}_0,\,\,(z,\eta) \in \mathscr{L}.$$
Note that for each $r \in (0,1]$ the deformed action of $\mathcal{H}_0$
is free on $\mathscr{L}$ and for each point of $\mathscr{L}_0$ there
is exactly one gauge orbit which goes through this point.
We now define $g_r$ for $r \in (0,1]$ to be the quotient metric
on $\mathscr{L}_0 \cong_r \mathscr{L} /_r 
\mathcal{H}_0$ of the $L^2$-metric $g$ on $\mathscr{L}$.
Here the $r$-parameter
takes account of the $r$-dependence of the action of $\mathcal{H}_0$
on $\mathscr{L}$. We denote by $\nabla_r \mathcal{A}_0$
the gradient of $\mathcal{A}_0$ with respect to the metric
$g_r$. To compute it,
note that t
he gradient of $\mathcal{A}^r$ with respect to the $L^2$-metric $g$ on
$\mathscr{L}$ for $r \in [0,1]$ is given by
$$\nabla \mathcal{A}^r (z,\eta)
=\left( \begin{array}{c}
i\partial_t z+(1-r)i X_{\bar{\eta}}(z)+r i X_\eta(z) \\
(1-r)\bar{\mu}(z)+r\mu(z)
\end{array}\right).$$
As in the proof of formula (\ref{grad2}) one shows
that for $(z,\eta) \in \mathscr{L}_0$ 
$$\nabla_r \mathcal{A}_0(z,\eta)=
(i\partial_t z+iL_z \eta+L_z \xi, \bar{\mu}(z)-\tau)$$
where $\xi \in \mathrm{Lie}(\mathcal{H}_0)$ is determined by
$$\partial_t \xi(t)=r^2\big(\mu(z(t))-\bar{\mu}(z)\big), 
\quad \int_0^1 \xi dt =0.$$
In particular,
$$\nabla_r \mathcal{A}_0(z,\eta)=r^2 \nabla_0 \mathcal{A}_0(z,\eta)
+(1-r^2)\nabla_1 \mathcal{A}_0(z,\eta), \quad 
r \in [0,1].$$
The energy of a flow line is defined to be
$$E(z,\eta)=\int_{-\infty}^\infty \int_0^1\big(||\partial_t z+L_z \eta||^2
+||\bar{\mu}(z)-\tau||^2\big)dt ds.$$
In order to prove Theorem C we first show the following lemma
which provides us with a uniform $L^\infty$-bound for $z_\nu$.
\begin{lemma}\label{Linfty}
There exists a compact set $K \subset \mathbb{C}^n$ such that for each flow
line $(z,\eta)$ of $\mathcal{A}^r_0$ for $r \in (0,1]$ 
whose energy is bounded, the trace of $z$ is contained in $K$, i.e. 
$$z(s,t) \in K, \quad s \in \mathbb{R},\,\, t \in S^1.$$
\end{lemma}
\textbf{Proof:} We first examine the Laplacian of the function
$|z|^2/2$. Observe that by (\ref{moment}) and Remark~\ref{proper}
for $z,z' \in \mathbb{C}^n$ 
$$\langle \mu(z),\mu(z')\rangle \geq 0.$$
Using this inequality and a computation similar to the one in the proof of 
\cite[Proposition 3.5]{cieliebak-gaio-salamon} we estimate
\begin{eqnarray}\label{cgs}
\Delta |z|^2/2&=&|\partial_s z+L_z \xi|^2+|\partial_t z+L_z \eta|^2\\
\nonumber
& &+2\big\langle \mu(z),
(r^2-1)\bar{\mu}(z)+r^2 \mu(z)-\tau \big\rangle\\ \nonumber
&\geq&2r^2\langle \mu(z),\mu(z)-\tau \rangle\\ \nonumber
&\geq&2r^2|\mu(z)|\big(|\mu(z)|-|\tau|\big).
\end{eqnarray}
We set
$$R:=\max_{|\mu(z)| \leq \tau}\{|z|\}, \quad
K:=\{z \in \mathbb{C}^n:|z| \leq 2R \}.$$
We assume by contradiction that there exists $(s_0,t_0) \in
\mathbb{R}\times S^1$ such that
$$z(s_0,t_0) \in \mathbb{C}^n \setminus K.$$
By the convexity property derived in (\ref{cgs}) it follows that either
for all $s \geq s_0$ or for all $s \leq s_0$ there exists 
$t_s$ with 
$$|z(s,t_s)| \geq |z(s_0,t_0)|.$$ 
We consider only the case where the above property holds for all $s \geq
s_0$, the treatment of the other case is completely analogous. Using the
assumption that the energy is bounded we conclude that there exist for
each $i \in \mathbb{N}$ real numbers $\sigma_i<\sigma^i$ with the
property that
$$\sigma_1 \geq s_0, \quad \sigma^i \leq \sigma_{i+1}\,\,
\forall \,\, i \in \mathbb{N}, \quad \sum_{i=1}^\infty (\sigma^i-\sigma_i)
=\infty$$
such that for each $s \in [\sigma_i,\sigma^i]$ for $i \in \mathbb{N}$ it 
exists $t'_s \in S^1$ such that
$$z(s,t'_s) \in K.$$
Abbreviating
$$\delta:=|z(s_0,t_0)|-R>0$$
we estimate
\begin{eqnarray*}
\infty &=&\sum_{i=1}^\infty \delta^2(\sigma^i-\sigma_i)\\
&\leq& \sum_{i=1}^\infty \int_{\sigma_i}^{\sigma^i}
\Bigg(\int_{t'_s}^{t_s}||\partial_t z+L_z \eta|| dt\Bigg)^2 ds\\
&\leq& \sum_{i=1}^\infty \int_{\sigma_i}^{\sigma^i}
\Bigg(\int_0^1||\partial_t z+L_z \eta|| dt\Bigg)^2 ds\\
&\leq& \sum_{i=1}^\infty \int_{\sigma_i}^{\sigma^i}
\int_0^1||\partial_t z+L_z \eta||^2 dtds\\
&\leq& \int_{-\infty}^\infty
\int_0^1||\partial_t z+L_z \eta||^2 dtds\\
&\leq&E(z,\eta)\\
&<&\infty.
\end{eqnarray*}
This contradiction proves the lemma. \hfill $\square$
\\ \\
\textbf{Proof of Theorem C: }By Theorem~\ref{embedding}
we may assume without loss of generality that $r_\nu>0$ for every
$\nu \in \mathbb{N}$. Choose $h_\nu \in H$ such that
$$0 \leq (h_\nu)_*(\eta_\nu)_i(0)< 2\pi, \quad 1 \leq i \leq k.$$
We will replace in the following $(z_\nu,\eta_\nu)$ by 
$(h_\nu)_*(z_\nu,\eta_\nu)$. By Lemma~\ref{Linfty} the maps $z_\nu$
are uniformly bounded in the $L^\infty$-topology.
The proof of Theorem C reduces now to elliptic
bootstrapping for the equations
\begin{eqnarray}\label{sver}
\partial_s z_\nu+i\partial_t z_\nu+L_{z_\nu} \xi_\nu+iL_{z_\nu} \eta_\nu
&=&0\\ \nonumber
\partial_s \eta_\nu+\bar{\mu}(z_\nu)=\tau\\ \nonumber
\partial_t \xi_\nu=r_\nu^2\big(\mu(z_\nu)-\bar{\mu}(z_\nu)\big)\\ \nonumber
\int_0^1 \xi_\nu dt=0.
\end{eqnarray}
Fix $N \in \mathbb{N}$. Using the fact that there exists a constant
$c_0=c_0(N)$ such that
$$||\partial_s z_\nu+i\partial_t z_\nu||_{L^2([-N,N]\times S^1)}
\leq c_0||z_\nu||_{W^{1,2}([-N-1,N+1]\times S^1)}$$
we deduce from the uniform $L^\infty$-bound on $z_\nu$ and the uniform
bound on $\eta_\nu(0)$ from (\ref{sver}) that there exists $c_1>0$ such
that
\begin{equation}\label{sob1}
||z_\nu||_{W^{1,2}([-N,N]\times S^1)}+||\eta_\nu||_{L^2([-N,N])}
+||\xi_\nu||_{L^2([-N,N]\times S^1)} \leq c_1.
\end{equation}
In order to get control over the higher derivatives it is useful to use
the original $\mathcal{H}_0$-action on $\mathscr{L}$ to put 
(\ref{sver}) into Coulomb gauge on the cylinder. Let
$\zeta_\nu \in C^\infty([-N,N]\times S^1,\mathfrak{t}^k)$ be a solution
of the following Neumann problem on the finite cylinder $[-N,N]\times S^1$
with mean value zero
$$\Delta \zeta_\nu(s,t)=r_\nu^2 d(\mu(z_\nu)-\bar{\mu}(z_\nu))\partial_s
z_\nu(s,t)=\partial_s \partial_t \xi_\nu
,\,\,s \in [-N,N],\,\,t \in S^1,$$
$$\partial_s \zeta_\nu(\pm N,t)=0,\,\,
t \in S^1,$$
$$\int_0^1 \zeta_\nu(s,t)dt=0,\,\, s \in [-N,N].$$
Define $g_\nu \in C^\infty([-N,N]\times S^1,T^k)$ as the map
from the interval $[-N,N]$ to the gauge group $\mathcal{H}_0$ which is
defined by the property
$$\partial_t (g_\nu^{-1}\partial_t g_\nu)=\zeta_\nu.$$
Define 
$(\tilde{z}_\nu,\tilde{\eta}_\nu,\tilde{\xi}_\nu) \in 
C^\infty([-N,N]\times S^1,\mathbb{C}\times \mathfrak{t}^k 
\times \mathfrak{t}^k)$ by
$$(\tilde{z}_\nu,\tilde{\eta}_\nu,\tilde{\xi}_\nu):=
(g_\nu)_*(z_\nu,\eta_\nu,\xi_\nu):=(\rho(g_\nu)z_\nu,\eta_\nu-
(g_\nu)^{-1}\partial_t g_\nu,\xi_\nu-(g_\nu)^{-1}\partial_s g_\nu).$$
Then $(\tilde{z}_\nu,\tilde{\eta}_\nu,\tilde{\xi}_\nu)$ is a solution of
the problem
\begin{eqnarray}\label{sver2}
\partial_s \tilde{z}_\nu+i\partial_t \tilde{z}_\nu
+L_{\tilde{z}_\nu}\tilde{\xi}_\nu+iL_{\tilde{z}_\nu}\tilde{\eta}_\nu=0\\
\nonumber
\partial_s \tilde{\eta}_\nu-\partial_t\tilde{\xi}_\nu
+(1-r_\nu^2)\bar{\mu}(z_\nu)+r_\nu^2\mu(z_\nu)=\tau\\ \nonumber
\partial_s \tilde{\xi}_\nu+\partial_t \tilde{\eta}_\nu=0\\ \nonumber
\int_0^1 \tilde{\xi}_\nu dt=0.
\end{eqnarray}
Combining (\ref{sob1}) with Lemma~\ref{lemma} below we see that there 
exists a constant $c_2>0$ such that
\begin{equation}\label{sob2}
||\tilde{z}_\nu||_{W^{1,2}([-N,N]\times S^1)}+
||\tilde{\eta}_\nu||_{L^2([-N,N]\times S^1)}
+||\tilde{\xi}_\nu||_{L^2([-N,N]\times S^1)} \leq c_2.
\end{equation}
By elliptic regularity for the Cauchy-Riemann operator it follows
from (\ref{sver2}) and (\ref{sob2}) that 
$(\tilde{z}_\nu,\tilde{\eta}_\nu,\tilde{\xi}_\nu)$ is uniformly bounded
on $[-N+1,N-1]\times S^1$ in the $W^{k,2}$-norm for every $k \in \mathbb{N}$. 
By the Sobolev embedding theorem and the theorem of Rellich it follows
that there exists a subsequence $\nu_j$ such that
$$(\tilde{z}_{\nu_j},\tilde{\eta}_{\nu_j},\tilde{\xi}_{\nu_j})
|_{[-N+1,N-1]\times S^1}
\stackrel{C^\infty}{\longrightarrow}_{j \to \infty}
(\tilde{z},\tilde{\eta},\tilde{\xi})$$
where $(\tilde{z},\tilde{\eta},\tilde{\xi})$ is a solution of
(\ref{sver2}) on $[-N+1,N-1]\times S^1$ for some $r \in [0,1]$.
Define $g \in C^\infty([-N+1,N-1]\times S^1,T^k)$ to be the map from
$[-N+1,N-1]$ to $\mathcal{H}_0$ defined by the property that
$g_*\tilde{\eta}$ is independent of the $t$-variable.
Then 
$$(z,\eta,\xi)=g_*(\tilde{z},\tilde{\eta},\tilde{\xi})$$ 
is a solution of (\ref{sver}). Moreover,
\begin{equation}\label{conv1}
(g \circ g_{\nu_j})_*(z_{\nu_j},\eta_{\nu_j},\xi_{\nu_j})
|_{[-N+1,N-1]\times S^1}
\stackrel{C^\infty}{\longrightarrow}_{j \to \infty}(z,\eta,\xi).
\end{equation}
Since $\eta_{\nu_j}$ for every $j \in \mathbb{N}$ and $\eta$ are
independent of the $t$-variable it follows that
$$\partial_t\big((g \circ g_{\nu_j})^{-1}\partial_t(g \circ g_{\nu_j}\big))
\stackrel{C^\infty}{\longrightarrow}_{j \to \infty}0.$$
It follows that
\begin{equation}\label{conv2}
g \circ g_{\nu_j} \stackrel{C^\infty}{\longrightarrow}_{j \to \infty}
\mathrm{id}.
\end{equation}
Combining (\ref{conv1}) and (\ref{conv2}) we conclude that
$$(z_{\nu_j},\eta_{\nu_j},\xi_{\nu_j})
|_{[-N+1,N-1]\times S^1}
\stackrel{C^\infty}{\longrightarrow}_{j \to \infty}(z,\eta,\xi).$$
Since $N \in \mathbb{N}$ was arbitrary the theorem follows. \hfill $\square$

\begin{lemma}\label{lemma}
There exists a constant $c>0$ with the following property.
Assume that $g_\pm \in C^\infty(S^1,\mathbb{R})$ and
$h \in C^\infty([-N,N] \times S^1,\mathbb{R})$ for $N \in \mathbb{N}$
satisfy 
$$\int_0^1 g_\pm(t)dt=0, \quad \int_0^1 h(s,t)dt=0,\,\, s \in [-N,N].$$
Suppose that $f \in C^\infty([-N,N]\times S^1, \mathbb{R})$ is a
solution of the Neumann problem on the finite cylinder
$[-N,N]\times S^1$ with mean value zero
$$\Delta f(s,t)=h(s,t),\,\,s \in [-N,N],\,\,t \in S^1,\quad 
\pm \partial_s f(\pm N,t)=g_\pm(t),\,\, t\in S^1,$$
$$\int_0^1 f(s,t)dt=0,\,\,
s \in [-N,N],$$
then
$$||f||_{W^{2,2}([-N,N]\times S^1)} \leq 
c\big(||h||_{L^2([-N,N]\times S^1)}+||g_+||_{W^{1,2}(S^1)}
+||g_-||_{W^{1,2}(S^1)}\big).$$
\end{lemma}
\textbf{Proof: }The proof is an exercise in partial
integration. We estimate
\begin{eqnarray}\label{partint}
& &\int_{-N}^N \int_0^1 \big((\partial^2_s f)^2+(\partial^2_t f)^2
+2(\partial_s \partial_t f)^2\big)dt ds \\ \nonumber
&=&\int_{-N}^N \int_0^1 \big((\partial^2_s f)^2+(\partial^2_t f)^2
-2(\partial^2_s \partial_t f)(\partial_t f)\big)dt ds\\ \nonumber
& &+\int_0^1\partial_t g_+(t)\partial_t f(N,t)dt+
\int_0^1\partial_t g_-(t)\partial_t f(-N,t)dt\\\nonumber
&\leq&\int_{-N}^N \int_0^1 \big((\partial^2_s f)^2+(\partial^2_t f)^2
+2(\partial^2_s f)(\partial^2_t f)\big)dt ds\\ \nonumber
& &+||\partial_t g_+||_2 ||\partial_t f(N,\cdot)||_2
+||\partial_t g_-||_2 ||\partial_t f(-N,\cdot)||_2\\ \nonumber
&\leq&||h||^2_2
+16||g_+||^2_{1,2}+16||g_-||^2_{1,2}+\frac{1}{16}||\partial_t f(N,\cdot)||^2_2
+\frac{1}{16}||\partial_t f(-N,\cdot)||^2_2.
\end{eqnarray}
Using the assumption that the mean value of $f(s,\cdot)$ vanishes for each 
$s \in [-N,N]$ we obtain the estimates
\begin{equation}\label{hoepo}
||f||_2 \leq ||\partial_t f||_2 \leq ||\partial^2_t f||_2,\quad
||\partial_s f||_2 \leq ||\partial_s \partial_t f||_2.
\end{equation}
Choose a smooth function $\beta \colon [-N,N] \to \mathbb{R}$ which
satisfies 
$$\beta(s)=0,\,\,s \leq N-1, \quad \beta(N)=1, \quad
0 \leq  \beta(s) \leq 1,\,\,s \in [-N,N],$$
$$\beta'(s) \leq 2,\,\, s\in [-N,N].$$
We estimate
\begin{eqnarray}\label{rand}
||\partial_t f(N,\cdot)||^2_2
&=&\int_{N-1}^N \partial_s \Bigg( \beta(s)
\int_0^1 (\partial_t f)^2(s,t)dt\Bigg)ds\\ \nonumber
&=&\int_{N-1}^N\int_0^1\beta'(s)(\partial_t f)^2(s,t)ds dt\\ \nonumber
& &+\int_{N-1}^N \int_0^1 2\beta(s) \partial_s \partial_t f(s,t)
\partial_t f(s,t)ds dt\\ \nonumber
&\leq&4||\partial_t f||_2^2+2||\partial_s \partial_t f||_2^2\\ \nonumber
&\leq&4||\partial_t^2 f||^2_2+2||\partial_s \partial_t f||_2^2.
\end{eqnarray}
For the last inequality we have used (\ref{hoepo}). Similarly one shows
\begin{equation}\label{rand2}
||\partial_t f(-N,\cdot)||^2_2 \leq 
4||\partial_t^2 f||^2_2+2||\partial_s \partial_t f||_2^2.
\end{equation}
Combining (\ref{partint}) with (\ref{rand}) and (\ref{rand2}) one
obtains
\begin{equation}\label{neumann}
||\partial_s^2 f||_2^2+||\partial_t^2 f||_2^2+
||\partial_s \partial_t f||_2^2 
\leq 2||h||_2^2+32||g_+||^2_{1,2}+32||g_-||^2_{1,2}.
\end{equation}
The lemma follows now by combining (\ref{neumann}) with (\ref{hoepo}).
\hfill $\square$

\end{document}